# Applied Trajectory Design for close-proximity operations of Asteroid CubeSat Missions

Claudio Bottiglieri*, Felice Piccolo†, Carmine Giordano‡ and Francesco Topputo§
*Politecnico di Milano, Milano, Italy, 20156*

**In this paper, a practical approach to the trajectory design for asteroid exploration missions with CubeSats is presented. When applied trajectories are sought, operative concerns and uncertainties affecting the spacecraft dynamics must be considered during the design process. Otherwise, trajectories that are possible on paper might become unfeasible when real-world constraints are considered. The risk of such eventualities leads to the urge of extending the trajectory design focus on the uncertainties affecting the dynamics and on the operative constraints derived by ground operations. This is especially true when targeting highly perturbed environments such as small bodies with low-cost solutions as CubeSats, whose capabilities in deep-space are still unknown. The case study presented is the Milani CubeSat which will be launched in 2024 with Hera in the frame of the AIDA mission.**

## Nomenclature

| | | |
|---|---|---|
| $a$ | = | semi-major axis, km |
| $e$ | = | eccentricity |
| $i$ | = | inclination, deg |
| $T$ | = | period, days |
| $M$ | = | asteroid mass, kg |
| $\mu$ | = | gravitational parameter, m$^3$/s$^2$ |
| $r$ | = | range, km |
| $\mathbf{r}_i$ | = | CubeSat position with respect to the i, km |
| $m$ | = | CubeSat mass, kg |
| $A$ | = | CubeSat area, m$^2$ |
| $C_r$ | = | reflectivity coefficient |
| $P_0$ | = | solar flux at 1 AU |

*Research Assistant, Dept. of Aerospace Science and Technology; claudio.bottiglieri@polimi.it.
†PhD Student, Dept. of Aerospace Science and Technology; felice.piccolo@polimi.it.
‡PostDoc Fellow, Dept. of Aerospace Science and Technology; carmine.giordano@polimi.it.
§Full Professor, Dept. of Aerospace Science and Technology; francesco.topputo@polimi.it. AIAA Senior Member.

| $\Delta$v | = | velocity change, m/s |
|---|---|---|
| $\Phi$ | = | State Transition Matrix |

Subscripts

| S | = | Sun |
|---|---|---|
| 1 | = | Didymos |
| 2 | = | Dimorphos |
| DS | = | Didymos–Sun |
| SRP | = | Solar Radiation Pressure |
| SE | = | Sun–Earth |
| r | = | position |
| v | = | velocity |

## I. Introduction

The first mission targeting a small body dates back to 1996 when NASA's Near Earth Asteroid Rendezvous (NEAR-Shoemaker) concluded the first non-incidental close rendezvous with asteroid (253) Mathilde. NEAR-Shoemaker also achieved the first soft-landing on (433) Eros' surface in 2001 [1]. From that moment, Close Proximity Operations (CPO) around small bodies have emerged as a relevant topic in space exploration. However, CPOs are non-trivial activities for spacecraft targeting small bodies. Indeed, the need for high accuracy in a risky environment, combined with the short time available for operations, makes CPOs challenging from both a conceptual and a technological point of view. A certain degree of autonomy or a strong effort from ground operators were the strategies used in the past to achieve a successful rendezvous. In Hayabusa 1 [2] and Hayabusa 2 [3], two recent successful missions to asteroids, sample operations were conducted autonomously by the probes, while critical operations for ESA's Rosetta mission required a timely dedication by the on-ground flight dynamics team [4]. However, an increasing number of space exploration missions will exploit CubeSats [5]: from MARCO [6] to the imminent CubeSat missions M-ARGO [7], NEA SCOUT [8] and the two Hera's CubeSats, Milani [9] and Juventas [10]. Thus, understanding the real capabilities of low-cost nanosatellites that operate close to a small body is of great interest.

Overcoming the technological challenges posed by CPOs is not always possible when low-cost solutions, such as CubeSats, are employed, and strategies developed for traditional spacecraft may not be suitable for small platforms. Until CubeSats autonomy becomes a reality, regular on-ground operations must be considered, together with the associated time delays. The time needed to perform orbit determination (OD), prepare the sets of instructions, validate them and upload them has to be considered while planning manoeuvres. Furthermore, having a full-time dedicated on-ground team is not possible because of the low-cost nature of CubeSats. In addition, the performances of miniaturized components (e.g., thruster magnitude and pointing accuracy) are often not comparable to standard ones, causing higher uncertainty. Consequently, trajectories that exist on paper might become unfeasible when real-world constraints are considered. Thus, the uncertainty on CubeSats capabilities and the highly perturbed environment close to small bodies show the need for a paradigm shift in the mission trajectory design.

In this work, a practical approach to trajectory design is presented. When considering highly perturbed environment and low-cost solutions as CubeSats, it must be clear that a navigation assessment to check for the flyability of the trajectories is essential. The aim of this work is to show how operative concerns and uncertainties affecting the dynamics must be considered during mission analyses and, if necessary, influence the trajectory design in a recursive fashion. This approach has been adopted in the case study reported in this paper, Hera's Milani CubeSat.

This paper is structured as follows: Section II introduces the Milani mission. The dynamical environment is detailed and the vehicle is presented. Section III shows the trajectory design constraints, strategy and results, while Section IV is focused on the navigation assessment of the designed trajectories. Finally, in Section V the work is concluded.

## II. The Milani Mission

The Milani mission, named after the mathematician Andrea Milani[*], is framed within the Hera mission of the European Space Agency (ESA). Hera is the ESA part of the Asteroid Impact & Deflection Assessment (AIDA) [11] international collaboration with NASA, which is responsible for the Double Asteroid Redirection Test (DART) [12] kinetic impactor spacecraft. Hera and DART have been conceived to be mutually independent. However, their value is increased when combined. The target of the Hera mission is the (65803) Didymos binary asteroid, which will be impacted by DART in September 2022. Following the impact, a crater is expected to appear on the secondary asteroid, Dimorphos. Hera will rendezvous with the asteroid in January 2027 carrying two CubeSats as opportunistic payloads (Juventas and Milani). Milani is scheduled to be released in March 2027 in proximity of the system. The main scientific goals of Milani are: 1) to obtain a detailed mapping and characterization of the asteroids, 2) to characterize the dust environment around them. The CubeSat will achieve these goals during the two main scientific phases: Far Range Phase (FRP) and Close Range Phase (CRP) [13]. At the end of CRP Milani will be injected into a Sun Synchronous Terminator Orbit (SSTO) and will attempt a landing on Dimorphos. This work focuses on the challenging design of the CRP.

### A. Dynamical Environment

*1. The Asteroid*

Didymos is a binary Near-Earth Asteroid (NEA) of S-type, discovered in 1996, formed by Didymos, or D1 (the primary) and Dimorphos, or D2 (the secondary). The reference model for the system is given in [14]. Table 1 and Table

---

[*]https://www.esa.int/Space_in_Member_States/Italy/Andrea_Milani_1948_2018

2 report the main data.

Table 1   Binary system orbital parameters (semi-major axis, eccentricity, inclination, period) [14].

| System Parameters | | | |
|---|---|---|---|
| a | e | i | T |
| 1.66446 AU | 0.3839 | 3.4083° | 770 days |

Table 2   Didymos and Dimorphos mass and spin periods properties [14].

| Asteroids Parameters | | | |
|---|---|---|---|
| $M_1$ | $M_2$ | $T_1$ | $T_2$ |
| 5.226×10$^{11}$ kg | 4.860×10$^9$ kg | 2.26 h | 11.92 h |

The orbital properties are retrieved from the kernels of the Hera mission[†]. In the reference model, Dimorphos and Didymos are assumed to share the same equatorial plane on which their relative motion occurs. Dimorphos is assumed to be in a tidally locked configuration with Didymos. The latter assumption can be relevant when observing some features on the secondary, such as the crater created by the DART impact, since it would be illuminated and visible only at certain geometries.

In this work, two reference frames have been used. A quasi-inertial frame called *DidymosEclipJ2000* which is centered in the system barycenter and has the *x-y* plane on the ecliptic at the epoch J2000 and the *z*-axis orthogonal to that plane. The axes are inertially fixed and the system can be considered inertial for intervals of time negligible with respect to Didymos heliocentric motion. The second one is a non-inertial reference frame in which trajectories are shown for clarity's sake. This frame has its *x-y* plane on D1's equator, with the *x*-axis aligned to the projection of the Sun vector on the equator and the *z*-axis aligned to the south pole of Didymos. This frame is also centered in the system barycenter and it is called *DidymosEquatorialSunSouth* (see Figure 1).

*2. Perturbations*

The dynamics of Milani in the proximity of the binary asteroid is modeled using the perturbed Restricted Four-Body Problem (R4BP), considering shape-based models for the asteroids only when it is relevant (break-even at 800 m [15]). The most important source of perturbation is the Sun, both because of its gravitational attraction and of Solar Radiation Pressure (SRP). The contribution of the Sun as a fourth body is modeled as

$$\mathbf{a}_{4body} = -\mu_S \left( \frac{\mathbf{r}_S}{r_S^3} - \frac{\mathbf{r}_{DS}}{r_{DS}^3} \right) \quad (1)$$

---
[†]https://www.cosmos.esa.int/web/spice/data Version 1.10 (Last retrieved on 08/11/21)

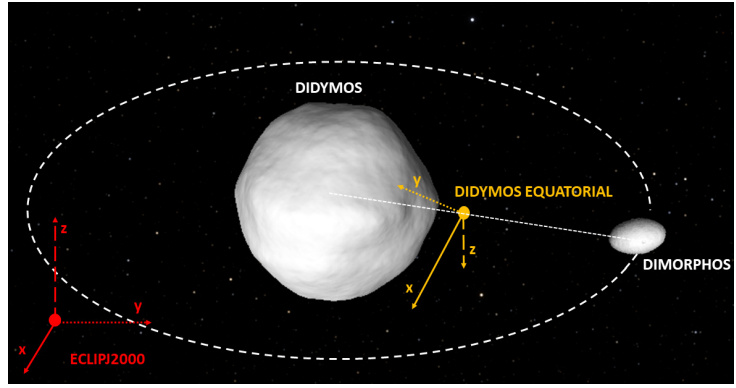

Fig. 1 Didymos geometry. The reference frames are highlighted. The red frame is the inertial Eclip2000 which corresponds to the quasi-inertial DidymosEclipJ2000 when centered in the system barycenter. The yellow frame is the DidymosEquatorialSunSouth.

where $\mu_S$ is the Sun gravity constant, $\mathbf{r}_S$ is the relative position of Milani with respect to the Sun and $\mathbf{r}_{DS}$ is the position of the system with respect to the Sun. Instead, the SRP has been modeled with a cannonball model[16].

$$\mathbf{a}_{SRP} = \frac{P_0}{c}\left(\frac{r_{SE}}{r_S}\right)^2 \frac{C_r A}{m} \frac{\mathbf{r}_S}{r_S} \qquad (2)$$

where, $P_0$ (1367 W/m$^2$) is the solar flux at 1 AU, $c$ is the speed of light, $r_{SE}$ is the Sun-Earth distance (1 AU), $C_r$ is the reflectivity coefficient of the CubeSat, $A$ is its equivalent surface area, and $m$ is its mass. The use of the cannonball model simplifies the mission design since it decouples the orbital dynamics from the spacecraft attitude.

*3. Equations of motion*

Considering the point-mass model for both asteroids, the equation of motion for Milani in the quasi-inertial DydimosEclipJ2000 reference frame can be written as

$$\ddot{\mathbf{r}} = -\mu_1 \frac{\mathbf{r}_1}{r_1^3} - \mu_2 \frac{\mathbf{r}_2}{r_2^3} + \mathbf{a}_{4body} + \mathbf{a}_{SRP} \qquad (3)$$

where $\mathbf{r}$ is the CubeSat position, $\mu_1$ and $\mu_2$ are the standard gravitational parameters of the primary and the secondary asteroid, while $\mathbf{r}_1$ and $\mathbf{r}_2$ are the CubeSat relative position with respect to the asteroids.

**B. The Vehicle**

Milani is a 6U CubeSat with 6 DOF maneuvering capabilities to control both translational and attitude motion. While attitude control is achieved autonomously on-board, orbital manoeuvres can be only triggered from ground. Table 3 reports the CubeSat parameters assumed for the SRP model.

**Table 3  Relevant parameters of the CubeSat (inputs to the Solar Radiation Pressure model)**

| Area [m$^2$] | C$_r$ | Mass [kg] |
|---|---|---|
| 0.51 | 1.25 | 12 |

*1. Instruments*

Milani's sensor suite is made of two Sun sensors, a star tracker, an IMU, a LIDAR and a navigation camera (NavCam). Additionally, as part of Milani's technological goals, the CubeSat has to exchange information and communicate with Hera via an inter-satellite link (ISL). This instrument also provides information on range and range-rate with respect to Hera. Thus, Milani optical navigation will also be supported by ISL measurements. Table 4 shows all the information assumed to simulate Milani navigation. NavCam uncertainties are expressed as function of the nominal range $r$. In the table, two sources of bias for the ISL measurements are considered: the typical instrument bias and the uncertainties on Hera's state.

**Table 4  Relevant assumed instruments parameter**

| Instrument | Noise on the measurement [1-sigma] | Bias |
|---|---|---|
| NavCam | $f_r = c_{r0} + c_{r1} r$  (range) | 2 m (range) |
|  | $f_\theta = \arctan\left(\frac{c_{\theta 0}}{r}\right) + c_{\theta 1} \sqrt{r}$  (angle) | 32.4 arcsec (angle) |
| ISL (Range) | 0.5 m | 150 m |
| ISL (Range-Rate) | 1.5 cm/s | 3 cm/s |

*2. The Scientific Payload: ASPECT*

Milani's platform accommodates three scientific payloads: ISL, VISTA [17], and ASPECT [18]. Among these, ASPECT is the most demanding from the trajectory design point of view. ASPECT is a passive payload, equipped with a visible to near-infrared hyperspectral imager, which will be used to perform global mapping of the asteroids with detailed observations of the DART crater on Dimorphos. Milani trajectories have been designed considering mainly ASPECT goals for both the scientific phases. In particular, during the close approach in CRP, Milani needs to image Dimorphos with a spatial resolution better than 1 m/pixel and observe the DART crater with a spatial resolution better than 0.5 m/pixel at phase angle (Sun–Asteroid–Milani angle) between [0–10] deg and [30–60] deg. The latter is the most challenging requirement and it is the main driver for the CRP design, while the mapping of whole the body will be a by-product of the crater observation arcs design. Indeed, the 0.5 m/pixel translates into a desired distance of 2780 m from D2. In addition, the requirement on the phase angle in between [0–10] deg, combined with the tidally locked nature of the system, reduces the admissible science region to portions of space 2 km far from the main asteroid with a consequent high risk of impacting the system.

# III. Trajectory Design

Milani trajectory design has been driven by the main scientific goals of the mission, but it has also been influenced by both technical and operational constraints. Differently from the technical constraints, which depend on the on-board instrument, Milani operational constraints are common to many space missions. Indeed, the challenges posed by those constraints can be considered ordinary effort for practical trajectory design studies, while in literature they are often ignored in preliminary analyses. Milani trajectories have been designed focusing on a realistic scenario, where maneuvering without the ground in the loop would be unfeasible, especially in proximity of harsh environments where uncertainties and perturbations play an important role. During the design phase, two solutions fulfilling all the requirements and constraints have been considered for Milani Close Range Phase. Their analysis are compared in Section IV.

## A. Payload and Ground Constraints

An important constraint on the design space is introduced by the passive nature of ASPECT and the NavCam. As a consequence, Milani must fly in the day-side of the system to perform both scientific observations and optical navigation. Therefore, the CubeSat needs to maneuver frequently to avoid going into the night-side, performing a repetition of hovering legs that span the day-side region back and forth. Although having a high maneuvering frequency can be costly, it might be desirable in response to highly nonlinear dynamical environments or to scientific needs. However, having daily or hourly maneuvering schedule is unrealistic from an operational point of view. In fact, the need to maneuver frequently is in contrast with the operational constraints associated with commanding the CubeSat from ground. The time needed to perform orbit determination and to prepare the correction maneuvers, and the working schedule of operative centers cannot ensure daily maneuvering, if not occasionally. Hence, to cope with these problems, an estimation of the Turn-Around-Time (TAT) has first to be made. The TAT is the time elapsed from the downlink of the observables to the uplink of the telecommands obtained after all on-ground operations are completed. The TAT can be considered the theoretical minimum time between two maneuvers commanded by ground. However, depending on the operative centers' working schedule, the real minimum time to fly a ballistic arc can be higher than the TAT. Furthermore, there could be an influence on the day of the week chosen for maneuvers. These details are particularly important when working on spacecraft formation missions around multiple small bodies, especially when phasing with other bodies is needed.

Milani OD will be performed simultaneously with Hera, for which a TAT of 48 hours has been considered. For this reason, no ballistic arcs with a duration of less than 48 hours can be designed to avoid open-loop maneuvering. Additionally, since Milani can communicate with ground only through Hera via the ISL, the CubeSat must use the same downlink/uplink window used by Hera, as no additional resources are foreseen for the CubeSat operations. Consequently, Milani's maneuvering schedule should be aligned with Hera's to exploit the information coming from ground as long as they are reliable. Thus, an additional constraint on the maneuvering schedule and frequency has to be considered since Milani should maneuver roughly every 3 and 4 days, close to Hera.

**B. Way Points and Key Points**

Previous work investigated trajectory options around Didymos for a CubeSat [15] and showed early-phase trajectories for Milani [19]. From these studies, the best design strategy for mid-range observations (10 km) emerged to be the *way-points strategy*. The way-point strategy is based on the selection of way points, which act as maneuver points between patched ballistic arcs. The geometry of the way points can be designed to fulfill the scientific requirements and meet the mission constraints as well, as seen with the Rosetta spacecraft about comet 67P/Churyumov-Gerasimenko [4]. A simple way-points strategy has been used to design Milani FRP, resulting in a repetition of hyperbolic orbits quasi-symmetric to the Sun direction. By placing the way points as far as possible, the spacecraft is forced to cover a larger distance in the same time (since the time of flight is fixed by the constraints), increasing its velocity and reducing its range from the target during the trajectory. Therefore, the relative way points placement can control the minimum range reached during the arc. However, this strategy is not suitable to fulfill the more challenging resolution requirements for CRP. Due to the constraints on the maneuvering frequency, which impose a minimum arc duration of 3 days, loop orbits quasi-symmetric to the Sun can reach a minimum distance with respect to Dimorphos of 7 km. The solution is to modify the way-points strategy with the introduction of *Key Points* [9]. A Key Point is the location at which the satellite can perform the desired scientific observation while fulfilling all the scientific requirements. Thus, the design focus moves from the maneuvers points geometry to the definition of the Key Points, one for each crater measurement. Figure 2 shows the Key Point design for the crater observations under the resolution and phase angle requirements described in Section II.B.2.

Consequently, the design of the maneuvers performed after the observations at the Key Points is done only to avoid going into the night-side. The resulting trajectory appears strongly asymmetrical in contrast with the previous scientific phase.

**C. Trajectory Solutions**

Two solutions have been found for Milani CRP. For convenience, they will be named "Option A" and "Option B" and their trajectories and geometry are shown in the remainder of this section.

*1. CRP-Option A*

Option A is an 8-point loop built with a 3-4 maneuvering pattern repeated four times. Figure 3 shows the projection of the loop orbit on the *x-y* and *y-z* planes in the *DidymosEquatorialSunSouth* reference frame. The left plots of Figure 3 show the time variation of phase angles and distances with respect to D1 and D2. Red horizontal dashed lines delimit the science admissible region for global mapping of D1 and D2, both in terms of phase angle (between 5 and 25 deg) and distance (between 2.780 and 4.572 km). A third dashed line is present in the bottom-right part of the figure representing the distance at which D1 saturates the NIR-FOV of ASPECT (1960 km). The time profiles show that the Key Point of each odd arc falls within the science admissible region for imaging the crater, both in terms of phase angle and distance.

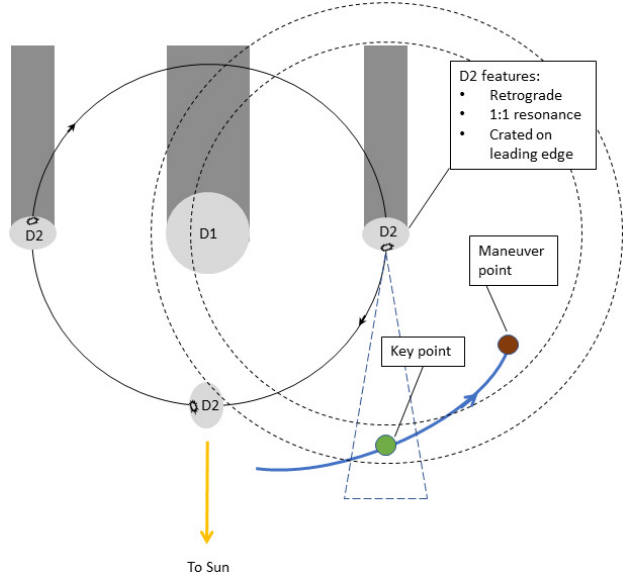

**Fig. 2  The choice of the "Key Point" where the CubeSat can image the crater on Dimorphos at optimal illumination conditions, at required distance and phase angle. The maneuver points of the observation arcs are chosen to avoid going into the night-side.**

*2. CRP-Option B*

Option B is a 6-point loop with a 7-4-7-7-3-7 days maneuvering pattern. The 7-days arcs can be seen as two arcs with a maneuvering frequency of 3 and 4 days, or vice versa, patched with a "zero magnitude" maneuver in between. Option B trajectory and geometry is shown in Figure 4. The main differences between the two options are the number of passage through the Key Points which is halved in Option B with respect to Option A, and the maneuvering frequency which is lower in Option B. Maneuvering each week causes a reduction of the number of arcs (and so of the number of possible scientific observations) to keep the duration of the phase constant. On the other hand, a greater maneuvering frequency can be extremely beneficial from the operational point of view, as it will be shown in Section IV.

## IV. Trajectory Analysis

Applied trajectory design goes past the determination of the nominal trajectory. Once in space, the spacecraft will not follow the nominal trajectory due to the uncertainties affecting the dynamics and control, which cannot be considered in the nominal design phase. Consequently, it is important to assess how much the real trajectory can differ from the nominal one and, in turn, to estimate the amount of navigation cost necessary to perform correction maneuvers.

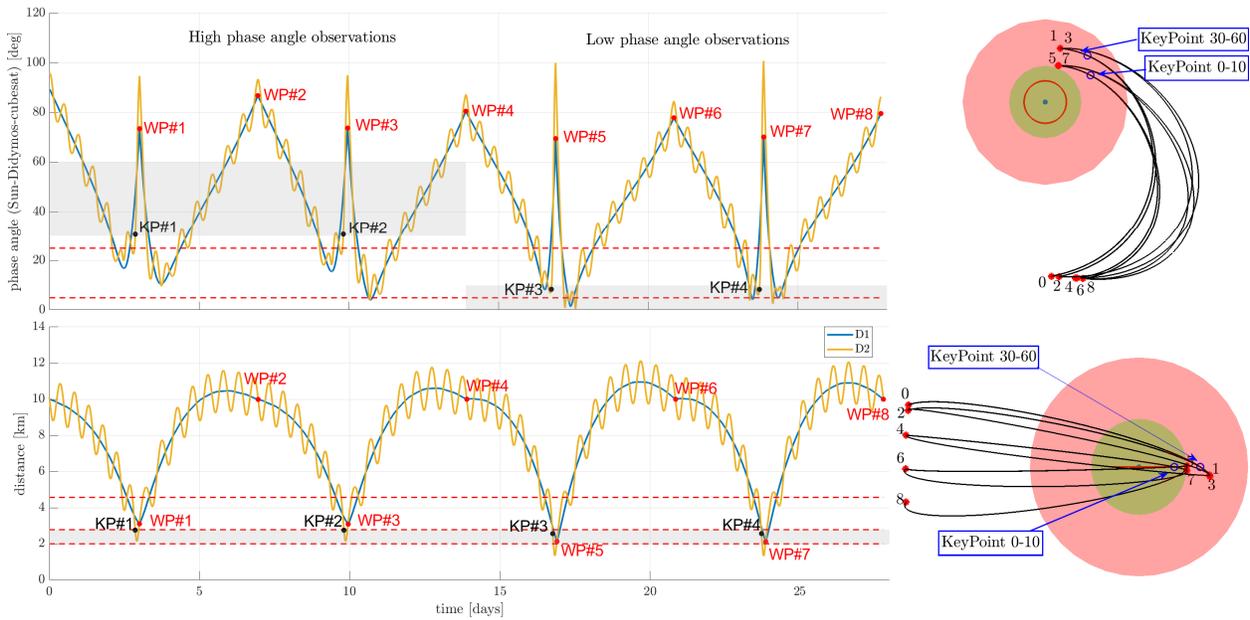

Fig. 3 **CRP-Option A orbit strategy: 8-point loop. Phase angle and distance time profile wrt D1 and D2 (upper and lower left). Key Points are depicted as black dots, Way Points are depicted as red dots; projection on DidymosEquatorialSunSouth *x*-*y* and *y*-*z* plane (upper and lower right). Grey areas are the requirements for the crater observations.**

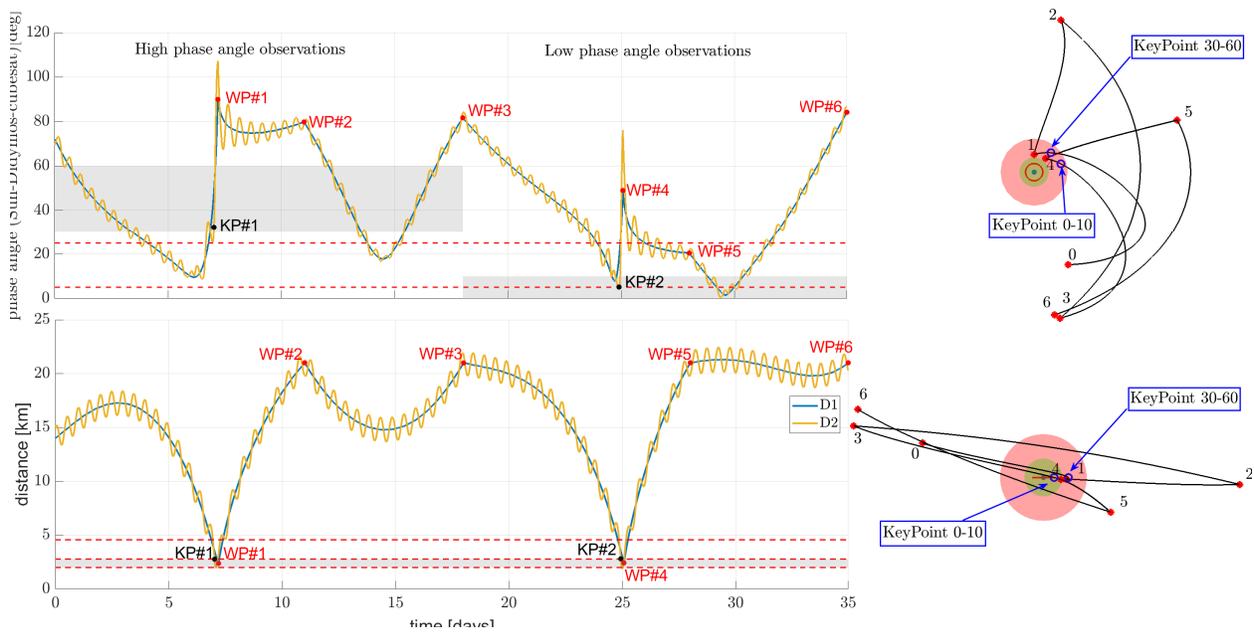

Fig. 4 **CRP-Option B orbit strategy: 6-point loop. Phase angle and distance time profile wrt D1 and D2 (upper and lower left). Key Points are depicted as black dots, Way Points are depicted as red dots; projection on DidymosEquatorialSunSouth *x*-*y* and *y*-*z* plane (upper and lower right). Grey areas are the requirements for the crater observations.**

This analysis measures the flyability of a trajectory which is of great interest for practical applications. Indeed, when flying in highly non-linear environments the difference between nominal and real trajectories can be not negligible, such that colliding with or escaping from the system can be a possibility. Thus, it is crucial to distinguish between nominal trajectories that exist only in the modeled reality and flyable orbits that are robust to uncertainties. The flyability of a trajectory can be assessed through two analyses: Knowledge Analysis and Dispersion Analysis.

### A. Knowledge Analysis

The knowledge analysis (KA) assesses the achievable level of accuracy in the spacecraft state knowledge. The term knowledge indicates the difference between the estimated trajectory and the real one. The level of accuracy of the state estimation is crucial since it is of interest in most of the mission operations. From the flight dynamic perspective, correction maneuvers are computed from the estimated deviation of the real trajectory with respect to the nominal one. Thus, the better the knowledge, the more effective the correction maneuvers. KA is based on the forward propagation of an initial knowledge represented as a Gaussian distribution centered in the nominal state $\mathbf{x}_0$ with a given initial covariance $\mathbf{P}_{K0}$. Uncertainties in the dynamics, like residual accelerations, thruster inaccuracies or uncertainties in the dynamical model, usually produce an increase of the knowledge covariance in time. The covariance can decrease only through an orbit determination process. During OD, pseudo-measurements are employed in a navigation filter to estimate the spacecraft state. The navigation of Milani is simulated considering radiometric measurements (range and range-rate between Hera and Milani) from the ISL and optical observations using landmarks based methods. Observables are processed by an Extended Kalman Filter (EKF) to estimate the state vector of the spacecraft. Biases are accounted for using Schimdt's formulation [20]. A dedicated analysis is performed for each mission phase. The filter is initialized considering the knowledge at the end of the previous phase, to which the additional uncertainty given by maneuvers is added. Both the state deviation and the covariance matrix are then propagated with the associated dynamics up to the first measurement epoch, where the estimated trajectory is updated. Propagation is linear and is based on the nominal state transition matrix (STM) as done in [21]. Proceeding in this way, the state estimates are sequentially updated as new measurements are processed, leading to the position and velocity knowledge profiles. This is repeated up to the final epoch of the orbit determination phase, obtaining the final spacecraft knowledge. To account for flight dynamics operations, a 48 hours TAT has been considered, during which the following operations are performed:

i) Data are sent from Milani to ground, through Hera.
ii) Measurements are processed and an OD solution is produced.
iii) Commands for the spacecraft are generated and validated on ground.
iv) Data are sent back to Milani, through Hera.

Furthermore, an additional margin of at least 1 hour is considered between data uplink and the execution of the maneuver. Therefore, the last observable that can be employed in the OD process is generated 49 hours before the following maneuver. This represents the so-called cut-off time (COT). However, letting the spacecraft fly for 49 hours without taking any measurement would increase the uncertainty unbearably, so the envisaged strategy is to take navigation images also during the COT. These images are then sent to Earth together with data taken during the following arc, in order to improve the post-maneuver knowledge by considering observables of the spacecraft state before the maneuver itself. Therefore, for each arc two covariances are propagated: one represents the actual ground knowledge for the current arc and considers only the measurement available before the COT, while the other represents the a-posteriori knowledge, after elaborating images taken during the COT and used for the following arc. The strategy is illustrated in Figure 5.

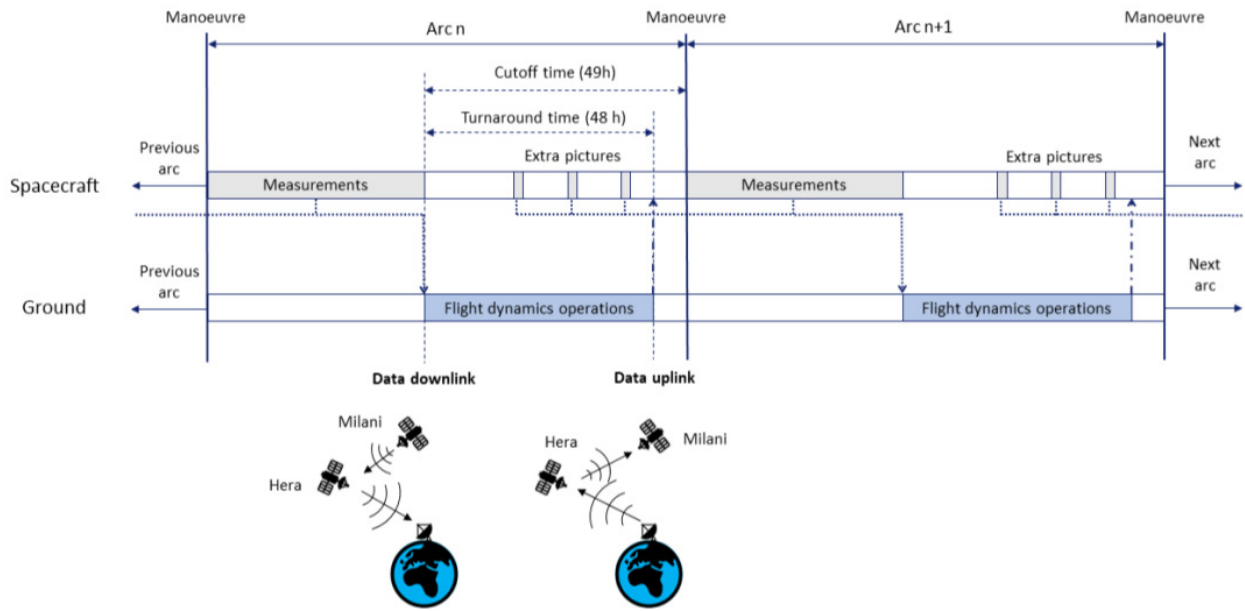

Fig. 5  Milani Knowledge Analysis timeline illustration.

The following assumptions are considered throughout the analysis:

- Orbit Determination is performed on the nominal trajectory.
- Navigation is considered to be relative to the Didymos System barycentre.
- A thrust error is considered with the uncertainty shown in Table 5.
- Range and range-rate measurements are supposed to be always available.
- SRP and residual acceleration are modeled as Gauss-Markov (GM) [21] processes, with the properties shown in Table 5.
- Uncertainty on the gravitational attraction is treated as a consider parameter. The assumed uncertainty at the beginning of the mission is of 3.5% of the nominal gravitational attraction of the system (35 $m^3/s^2$). This value was taken from Horizon Didymos reference model, properly reduced considering the contribution that DART and Hera will give for the gravity field estimation before Milani release. CRP takes places a month after the release when the uncertainty over the gravity is reduced to $10^{-4}$ m$^3$/s$^2$.

- For CRP-Option A the measurements plan is the following.
    - ISL measurement always start 6h after the beginning of the arc and end 3h before the cut-off-time. The frequency for range measurements is of 3h while for range-rate 1h.
    - 7 optical measurements are always taken in a single arc. 4 of them before the cut-off-time for the OD of the current arc, the other 3 after the cut-off-time and they will be used for the OD of the next arc.
- For CRP-Option B the measurements plan is the following.
    - ISL measurement always start 6h after the beginning of the arc and end 3h before the end of the arc. The frequency for range measurements is of 3h while for range-rate 1h.
    - 12 optical measurements are always taken in a single arc. 7 of them before the cut-off-time for the OD of the current arc, the other 5 after the cut-off-time and they will be used for the OD of the next arc.

Table 5   Uncertainties on the dynamics

| Uncertainty | Value [1-sigma] |
| --- | --- |
| Gravity uncertainty | $10^{-4}$ m$^3$/s$^2$ |
| Thruster | 1.67 % (magnitude),   0.67 deg (direction) |
| SRP | 8 % (magnitude),   1 day (correlation time) |
| Residual acceleration | $5 \times 10^{-9}$ m/s$^2$,   1 day (correlation time) |

A different measurement plan has been conceived for the two design options. The knowledge does not improve significantly by increasing the number of measurements further.

**B. Dispersion Analysis**

The Dispersion Analysis (DA) shows the statistics of the nominal trajectory from the real one. Indeed, the term *dispersion* indicates the stochastic distance between the two. DA can be performed with a Monte Carlo analysis where all the uncertainties affecting the dynamics are considered. For each run, the initial nominal state $\mathbf{x}_0$ is perturbed with an initial uncertainty represented as a Gaussian distribution centered in the nominal state with a given initial covariance $\mathbf{P}_{D0}$. This initial perturbed state is propagated under the effect of the perturbed dynamic up to the COT when the observables are downlinked. Here the estimated deviation is simulated by perturbing the CubeSat state with the knowledge and then propagated up to the correction time to get the estimated deviation needed for the correction maneuver evaluation. Instead, the real trajectory is propagated up to the correction time when it will receive the correction maneuver before proceeding with the next leg. Usually, for deep-space operations a correction maneuver is scheduled each week. For CPOs, a correction is usually needed at the earliest convenience. Therefore correction maneuvers can be scheduled concurrently with nominal maneuvers. As OD improves knowledge, correction maneuvers reduce dispersion. In order to estimate them a *guidance flow* is needed. Figure 6 shows the approach followed.

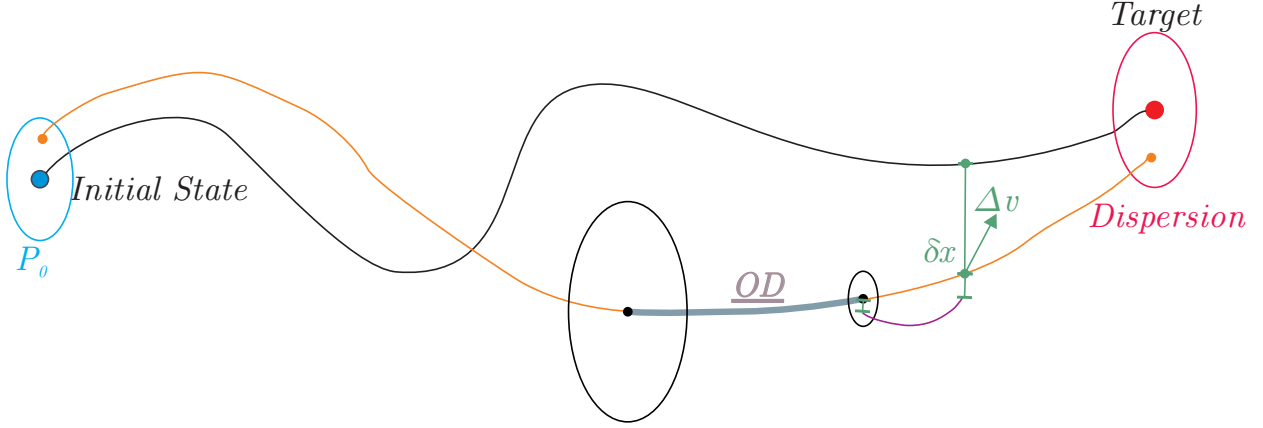

**Fig. 6   Dispersion analysis approach. The black line is the nominal trajectory, the orange line is one of the possible real trajectories and the purple one is the relative estimated trajectory (figure taken from [22]).**

In Milani the differential guidance, a commonly used guidance strategy for interplanetary missions [23], is applied. In this formulation, the whole trajectory is subdivided in different legs. At the extremal points of a single leg, two maneuvers are applied to cancel both the position and the velocity deviations on the final leg point. However, the final impulse is usually not applied in practice, since at the time of arrival at the final point a new maneuver is calculated in a receding horizon approach. Based on this, the maneuver can be computed by minimizing the deviations from the nominal state at the final point in a least square residual sense. Thus, the maneuver that must be applied at the time $t_j$ in order to cancel out the deviations at time $t_{j+1}$ is computed as:

$$\Delta \mathbf{v}_j = - \left( \mathbf{\Phi}_{rv}^T \mathbf{\Phi}_{rv} + q \mathbf{\Phi}_{vv}^T \mathbf{\Phi}_{vv} \right)^{-1} \left( \mathbf{\Phi}_{rv}^T \mathbf{\Phi}_{rr} + q \mathbf{\Phi}_{vv}^T \mathbf{\Phi}_{vr} \right) \widetilde{\delta \mathbf{r}}_j - \widetilde{\delta \mathbf{v}}_j \qquad (4)$$

where $\delta \mathbf{r}_j$ and $\delta \mathbf{v}_j$ are the estimated position and velocity deviation, $\mathbf{\Phi}_{rr}, \mathbf{\Phi}_{rv}, \mathbf{\Phi}_{vr}, \mathbf{\Phi}_{vv}$ are the 3-by-3 blocks of the STM $\mathbf{\Phi}$ from time $t_j$ to time $t_{j+1}$ associated to the nominal trajectory, and $q$ adjusts the dimensions. Control impulses are applied at nominal maneuvering points only, by superimposing the correction maneuver to the nominal ones. This procedure is repeated up to the final time. The estimation of the total navigation cost is obtained as the sum of all the maneuvers. The assumptions made during the DA are summarized and listed below:

- Each arc is initialized considering the dispersion at the end of the previous phase, to which the additional uncertainty given by maneuvers is added.
- A thrust misalignment of 5% in magnitude and 2 deg in angle (3-sigma) is considered.
- SRP and residual acceleration are modeled as Guass-Markov process with the same properties shown in Table 5.
- Uncertainty on the gravitational attraction is considered with the same properties shown in Table 5.

### C. Results

Knowledge and dispersion analysis results for both the CRP options are shown in this section. Figure 7 and Figure 9 show the position knowledge for Option A and Option B, respectively, while Figure 8 and Figure 10 show the velocity knowledge. For both options, the knowledge computed on-ground at the maneuver instants is in the order of a few hundreds of meters for the position and less than a cm/s for the velocity. No particular differences can be highlighted between the two options. Contrarily, the two solutions differ greatly from the dispersion point of view, as shown in Figure 11 and Figure 12. The figures show both the absolute dispersion (black) and relative dispersion in percentage (red) as the ratio between the absolute dispersion and the distance between the CubeSat and the asteroids. Figure 11 clearly shows the unflyability of the CRP-Option A. The error introduced by the starting maneuver is so high that the CubeSat reaches a relative dispersion of 33% (1-sigma) before any possibility of correcting such error. A relative dispersion of 33% (1-sigma) reveals that there is a non-negligible collision probability. It is evident how a correction maneuver is needed before approaching the system at few km to reduce the relative dispersion. Such precision can be obtained only with arcs of at least 7 days that allow for a correction maneuver in the middle when no concurrent expensive, and so prone to relevant uncertainties, maneuvers are scheduled. The navigation costs for this option are not shown since they are unreliable. Contrarily, although CRP-Option B dispersion is high, it shows how this trajectory can be considered reliable at least in terms of collision risk. For this reason, in absence of not-negligible collision risk, the figure shows the dispersion analysis for the entire mission phase.

The results of the navigation assessment clearly shows that Option A is unfeasible, even if its nominal trajectory exists according to the dynamical model considered. In absence of these results, the mission analysis would have presented a trajectory unflyable during operations. Although Option B is still risky, its design has been conceived to be as robust as possible to uncertainties to guarantee the flyability of the trajectory and so it is a feasible option for Milani.

## V. Conclusions

In this work, a practical approach to the trajectory design is presented. The case study is the Close Proximity Operation of Milani, a CubeSat mission around a binary asteroid. Milani mission analysis shows how operational constraints can be relevant to design flyable trajectories. Furthermore, operating in strongly perturbed environments, like about small bodies, leads to the urge of extending the ordinary trajectory design focus to reliability analyses on the designed nominal orbits. Milani dispersion analysis shows how existing nominal orbits can result unflyable when considering uncertainties on the maneuver execution or on the dynamical model considered.

This work highlights the need for a paradigm shift on the ordinary trajectory design strategy to a more practical approach. In order to be more aligned to the real world, operational constraints need to be considered and the impact of uncertainties affecting the dynamics should be assessed especially in highly perturbed environment such as small bodies. A navigation assessment should be performed to check for the flyability of the trajectory that can exists on paper, but

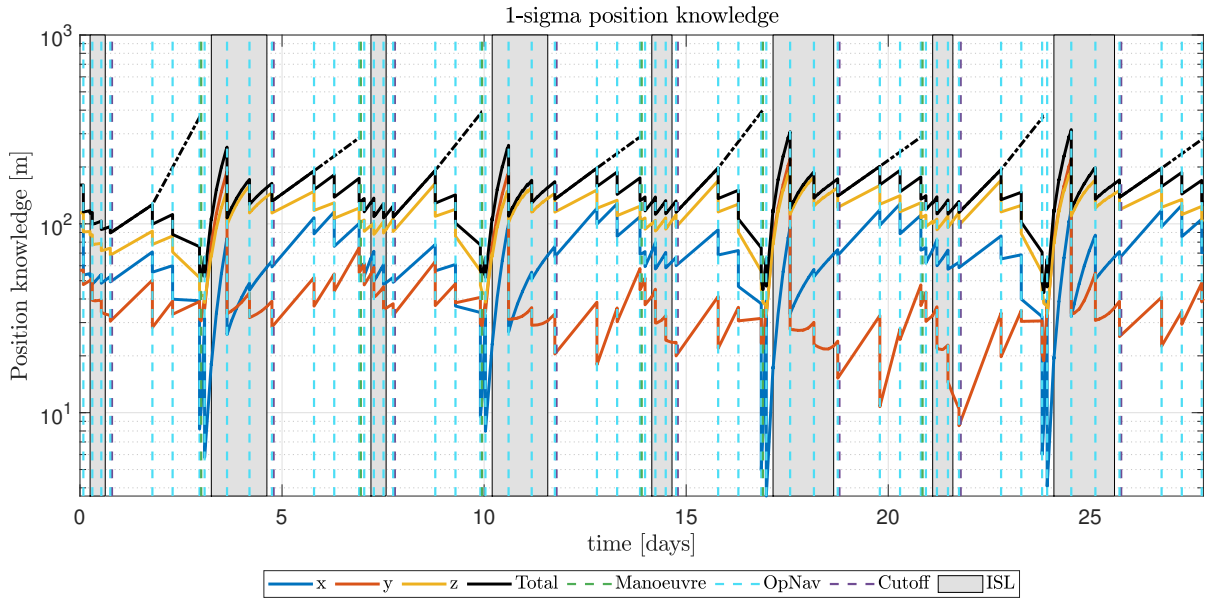

**Fig. 7  1-sigma position knowledge during the CRP-Option A. The dashed black line is the knowledge obtained considering only observables generated before COT, i.e., the knowledge known from ground.**

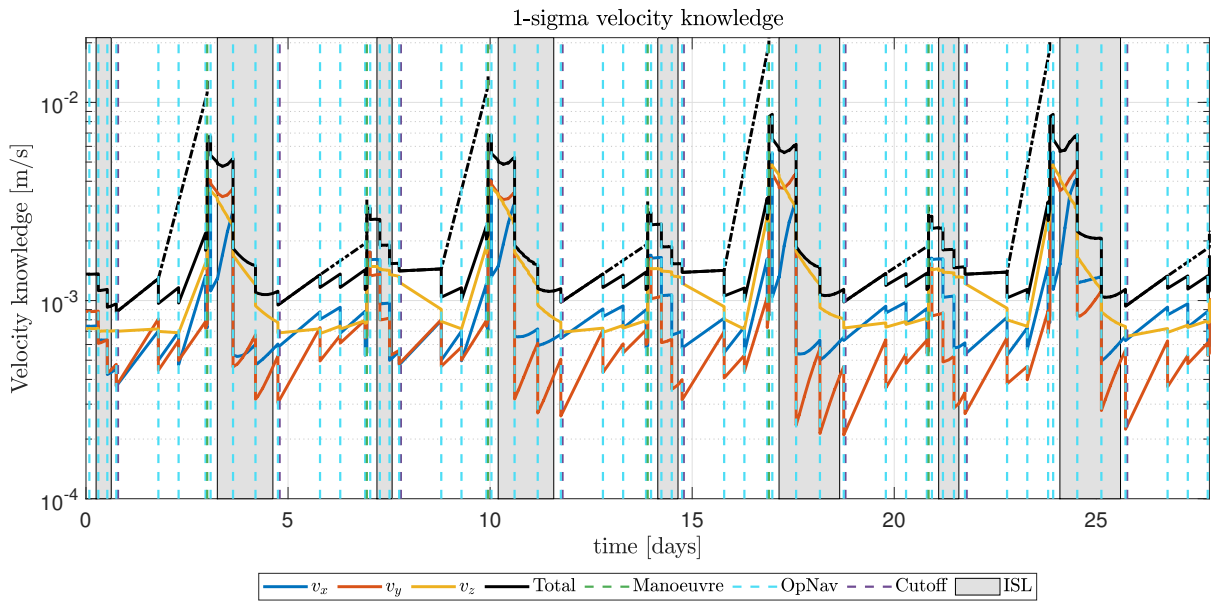

**Fig. 8  1-sigma velocity knowledge during the CRP-Option A. The dashed black line is the knowledge obtained considering only observables generated before COT, i.e., the knowledge known from ground.**

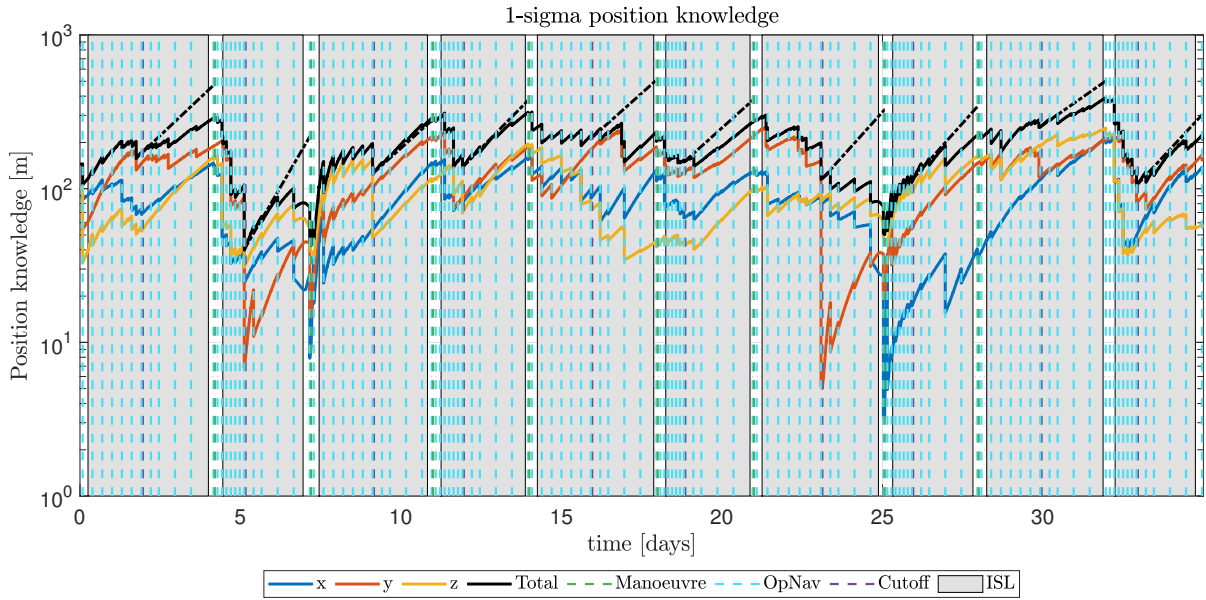

**Fig. 9  1-sigma position knowledge during the CRP-Option B. The dashed black line is the knowledge obtained considering only observables generated before COT, i.e., the knowledge known from ground.**

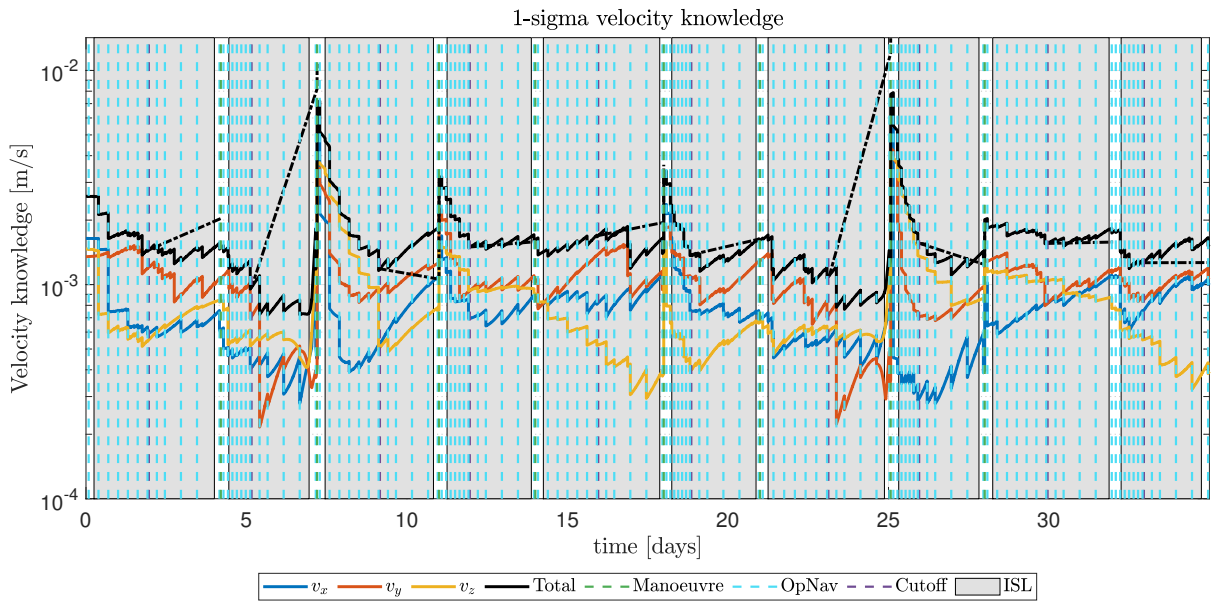

**Fig. 10  1-sigma velocity knowledge during the CRP-Option B. The dashed black line is the knowledge obtained considering only observables generated before COT, i.e., the knowledge known from ground.**

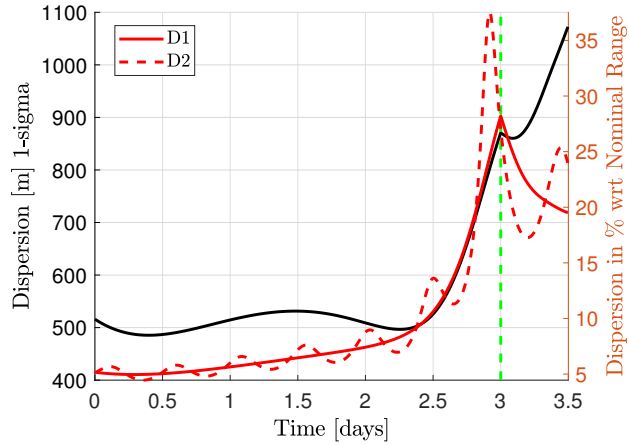

**Fig. 11 Dispersion Analysis results for CRP-Option A [1-sigma]. Absolute dispersion on the left (black), % value with respect to the nominal range on the right (red). Dashed green lines are the maneuvers. Maneuvers are shown as dashed green lines. The graph stops after the first maneuver since a non-negligible collision risk is reached at that point.**

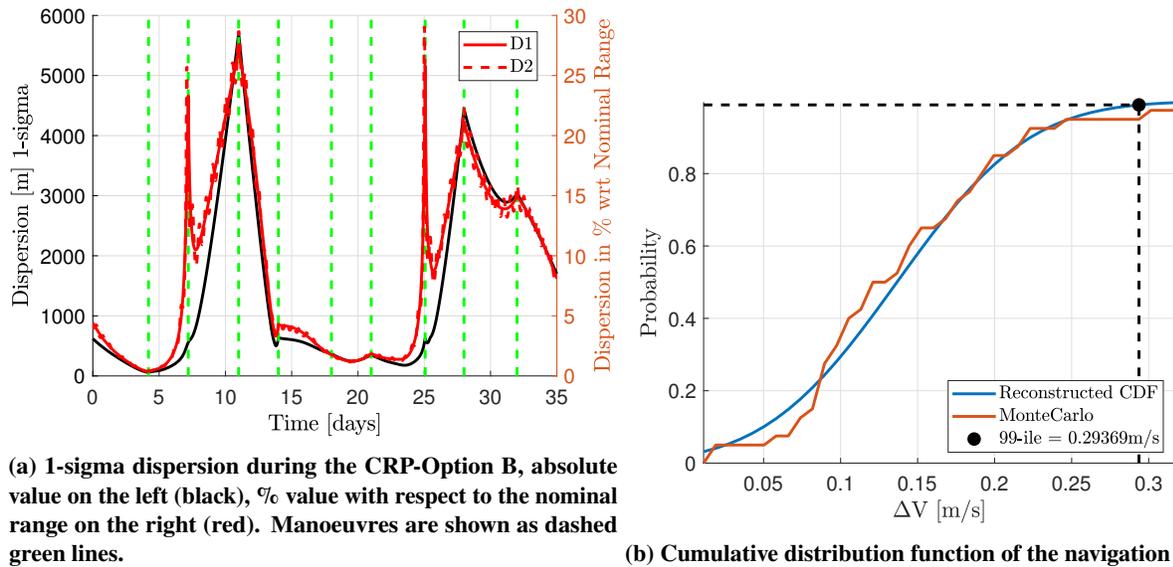

**(a) 1-sigma dispersion during the CRP-Option B, absolute value on the left (black), % value with respect to the nominal range on the right (red). Manoeuvres are shown as dashed green lines.**

**(b) Cumulative distribution function of the navigation costs.**

**Fig. 12 Dispersion Analysis results for CRP-Option B**

can be highly affected by uncertainties and thus be unflyable. Therefore, trajectories should be reiterated in case their flyability cannot be ensured against uncertainties.

## Acknowledgment

Part of the work presented in this paper has been performed under Tyvak International contract within ESA Contract No. 1222343567/62/NL/GLC. The authors would like to acknowledge the support received by the whole Milani Consortium.

## References

[1] Veverka, J., Farquhar, B., Robinson, M., Thomas, P., Murchie, S., Harch, A., Antreasian, P., Chesley, S., Miller, J., Owen, W., et al., "The landing of the NEAR-Shoemaker spacecraft on asteroid 433 Eros," *Nature*, Vol. 413, No. 6854, 2001, pp. 390–393. https://doi.org/10.1016/S0094-5765(02)00098-X.

[2] Furfaro, R., Cersosimo, D., and Wibben, D. R., "Asteroid Precision Landing via Multiple Sliding Surfaces Guidance Techniques," *Journal of Guidance, Control, and Dynamics*, Vol. 36, No. 4, 2013, pp. 1075–1092. https://doi.org/10.2514/1.58246.

[3] Yoshikawa, K., Sawada, H., Kikuchi, S., Ogawa, N., Mimasu, Y., Ono, G., Takei, Y., Terui, F., Saiki, T., Yasuda, S., et al., "Modeling and analysis of Hayabusa2 touchdown," *Astrodynamics*, Vol. 4, No. 2, 2020, pp. 119–135. https://doi.org/10.1007/s42064-020-0073-x.

[4] Accomazzo, A., Ferri, P., Lodiot, S., Pellon-Bailon, J.-L., Hubault, A., Porta, R., Urbanek, J., Kay, R., Eiblmaier, M., and Francisco, T., "Rosetta operations at the comet," *Acta Astronautica*, Vol. 115, 2015, pp. 434–441. https://doi.org/10.1016/j.actaastro.2015.06.009.

[5] Walker, R., Binns, D., Bramanti, C., Casasco, M., Concari, P., Izzo, D., Feili, D., Fernandez, P., Fernandez, J. G., Hager, P., et al., "Deep-space CubeSats: thinking inside the box," *Astronomy & Geophysics*, Vol. 59, No. 5, 2018, pp. 24–30. https://doi.org/10.1093/astrogeo/aty232.

[6] Martin-Mur, T. J., and Young, B., "Navigating MarCO, the first interplanetary CubeSats," *International Symposium on Space Flight Dynamics*, Melbourne, Australia, 2019.

[7] Topputo, F., Wang, Y., Giordano, C., Franzese, V., Goldberg, H., Perez-Lissi, F., and Walker, R., "Envelop of reachable asteroids by M-ARGO CubeSat," *Advances in Space Research*, Vol. 67, No. 12, 2021, pp. 4193–4221. https://doi.org/10.1016/j.asr.2021.02.031.

[8] McNutt, L., Johnson, L., Kahn, P., Castillo-Rogez, J., and Frick, A., "Near-Earth asteroid (NEA) scout," *AIAA Space 2014 Conference and Exposition*, AIAA Paper 2014–4435, San Diego, CA, 2014.

[9] Bottiglieri, C., Piccolo, F., Rizza, A., Giordano, C., Pugliatti, M., Franzese, V., Ferrari, F., and Topputo, F., "Trajectory design

and orbit determination for Hera's Milani CubeSat," *AAS/AIAA Astrodynamics Specialist Conference*, AAS 21–667, Big Sky, CA, 2021.

[10] Karatekin, O., Goldberg, H., Prioroc, C.-L., and Villa, V., "Juventas: Exploration of a binary asteroid system with a CubeSat," International Astronautical Federation Paper IAC-19.A3.4B.6.x54804, Washington DC, 2019.

[11] Cheng, A., Atchison, J., Kantsiper, B., Rivkin, A., Stickle, A., Reed, C., Galvez, A., Carnelli, I., Michel, P., and Ulamec, S., "Asteroid Impact and Deflection Assessment mission," *Acta Astronautica*, Vol. 115, 2015, pp. 262–269. https://doi.org/10.1016/j.actaastro.2015.05.021.

[12] Sarli, B. V., Atchison, J. A., Ozimek, M. T., Englander, J. A., and Barbee, B. W., "Double Asteroid Redirection Test Mission: Heliocentric Phase Trajectory Analysis," *Journal of Spacecraft and Rockets*, Vol. 56, No. 2, 2019, pp. 546–558. https://doi.org/10.2514/1.A34108.

[13] Bottiglieri, C., Piccolo, F., Rizza, A., Pugliatti, M., Franzese, V., Giordano, C., Ferrari, F., and Topputo, F., "Mission Analysis and Navigation Assessment for Hera's Milani CubeSat," *4S Symposium*, Vilamoura, Portugal, 2022.

[14] Naidu, S., Benner, L., Brozovic, M., Nolan, M., Ostro, S., Margot, J., Giorgini, J., Hirabayashi, T., Scheeres, D., Pravec, P., Scheirich, P., Magri, C., and Jao, J., "Radar observations and a physical model of binary near-Earth asteroid 65803 Didymos, target of the DART mission," *Icarus*, Vol. 348, 2020, p. 113777. https://doi.org/10.1016/j.icarus.2020.113777.

[15] Ferrari, F., Franzese, V., Pugliatti, M., Giordano, C., and Topputo, F., "Trajectory Options for Hera's Milani CubeSat Around (65803) Didymos," *The Journal of the Astronautical Sciences*, Vol. 68, No. 4, 2021, pp. 973–994. https://doi.org/10.1007/s40295-021-00282-z.

[16] Utashima, M., "Spacecraft orbits around asteroids for global mapping," *Journal of Spacecraft and Rockets*, Vol. 34, No. 2, 1997, pp. 226–232. https://doi.org/10.2514/2.3197.

[17] Palomba, E., Longobardo, A., Dirri, F., Zampetti, E., Biondi, D., Saggin, B., Bearzotti, A., and Macagnano, A., "VISTA: A $\mu$-Thermogravimeter for Investigation of Volatile Compounds in Planetary Environments," *Origins of Life and Evolution of Biospheres*, Vol. 46, No. 2, 2016, pp. 273–281. https://doi.org/10.1007/s11084-015-9473-y.

[18] Kohout, T., Näsilä, A., Tikka, T., Granvik, M., Kestilä, A., Penttilä, A., Kuhno, J., Muinonen, K., Viherkanto, K., and Kallio, E., "Feasibility of asteroid exploration using CubeSats—ASPECT case study," *Advances in Space Research*, Vol. 62, No. 8, 2018, pp. 2239–2244. https://doi.org/10.1016/j.asr.2017.07.036.

[19] Ferrari, F., Franzese, V., Pugliatti, M., Giordano, C., and Topputo, F., "Preliminary mission profile of Hera's Milani CubeSat," *Advances in Space Research*, Vol. 67, No. 6, 2021, pp. 2010–2029. https://doi.org/10.1016/j.asr.2020.12.034.

[20] Simon, D., *Optimal State Estimation: Kalman, H Infinity, And Nonlinear Approaches*, Wiley, 2006, pp. 309–312.

[21] Tapley, B. D., Schutz, B. E., and Born, G. H., *Statistical Orbit Determination*, Elsevier, 2004, Chap. 4.

[22] Giordano, C., "Analysis, Design, and Optimization of Robust Trjaectories for Limited-Capability Small Satellites," Ph.D. thesis, Politecnico di Milano, 2021. URL http://hdl.handle.net/10589/177695.

[23] Dei Tos, D. A., Rasotto, M., Renk, F., and Topputo, F., "LISA Pathfinder mission extension: A feasibility analysis," *Advances in Space Research*, Vol. 63, No. 12, 2019, pp. 3863–3883. https://doi.org/10.1016/j.asr.2019.02.035.